\documentclass[12pt]{article}

\usepackage{epsfig}
\usepackage{amssymb}

\newcommand{\VS}{\vspace{.5cm}}
\newcommand{\Vs}{\vspace{.3cm}}
\newcommand{\vs}{\vspace{.15cm}}

\newcommand{\mod}{\mbox{mod }}

\newcommand{\ZZ}{{\mathbb{Z}}}
\newcommand{\QQ}{{\mathbb{Q}}}

\title{Thue equations and torsion groups of elliptic curves}
\author{Irene Garc\'{\i}a--Selfa$^*$ \and Jos\'e M. Tornero\footnote{Both authors supported by FQM 218 (JdA), MTM2007-66929 and FEDER (MEC).}}
\date{June, 2007}

\begin{document}

\maketitle

\begin{abstract}
A new characterization of rational torsion subgroups of elliptic curves is found, for points of order greater than $4$, through the existence of solution for systems of Thue equations.
\end{abstract}

\VS



\section{Introduction}

In this paper we consider elliptic curves defined over $\QQ$. As it is known \cite{Cassels, Silverman}, each one of such curves is birationally equivalent to one, say $E$; given by an equation of the type
$$
E: Y^2=X^3+AX+B, \quad \mbox{ with }A,B, \in \ZZ;
$$
called short Weierstrass form, where it must hold $\Delta = 16(4A^4+27B^2) \neq 0$. The set of rational points of its projective clausure, noted $E(\QQ)$, is a finitely generated abelian group (Mordell--Weil Theorem \cite{Mordell, Weil}) and its torsion part, noted $T(E(\QQ))$ has been exhaustively described by Mazur \cite{Mazur, MazurIHES} as isomorphic to one of the following groups:
$$
\begin{array}{lcl}
\ZZ / n \ZZ & \mbox{ for } & n=1,2,...,10,12; \\
\ZZ / 2 \ZZ \times \ZZ / 2n \ZZ & \mbox{ for } & n=1,2,3,4.
\end{array}
$$

In a previous paper of ours \cite{GT} we proved the following result, charactering the non--trivial torsion subgroups by means of the (non--)existence of solution for a system of diophantine equations:

\vs

\noindent {\bf Theorem.--} Let $E:Y^2=X^3+AX+B$ be an elliptic curve, with $A,B \in \ZZ$. For every $n \in \{2,3,4,5,7,8,9\}$ there are, at most, $4$ quasi--homogeneous polynomials $P_n,Q_n,R_n,S_n \in \ZZ [z_1,...,z_4]$ such that $E$ has a rational point of order $n$ if and only if there exists an integral solution for the system
$$
\Sigma_n : \left\{ 
\begin{array}{rcl} 
P_n(z_1,...,z_4) &=& 6^2 \cdot A\\ 
Q_n(z_1,...,z_4) &=& 6^3 \cdot B \\
R_n(z_1,...,z_4) &=& 0 \\
S_n(z_1,...,z_4) &=& 0 \\
\end{array} \right.
$$

\vs

\noindent {\bf Remark.--} Only orders which are prime or pure prime powers are considered
because the remaining cases may be solved by joining systems, according to the factorization of the order we are interested on.

\vs

\noindent {\bf Remark.--} Here we show some more precise information regarding the systems $\Sigma_n$

\begin{tabular}{llc}
\\
{\bf Case} & {\bf Equations and variables} & {\bf Max. degree} \\
\hline \\
$n=2$ & $2$ equations in $2$ variables & $2$ and $3$ \\
$n=3$ & $2$ equations in $2$ variables & $6$ \\
$n=4$ & $2$ equations in $2$ variables & $3$ \\
$n=5$ & $3$ equations in $3$ variables & $4$ \\
$n=7$ & $4$ equations in $4$ variables & $4$ \\
$n=8$ & $3$ equations in $3$ variables & $6$ \\
$n=9$ & $3$ equations in $3$ variables & $9$\\
\\
\end{tabular}

The case $n=2$ has two different maximum degrees because, in fact, two different systems have to be used: one for detecting two--torsion points, and another for knowing whether there are one or three such points.

\vs

The systems $\Sigma_n$, for $n>4$ were to some extent unsatisfactory becuase they were quite heterogeneous: all systems were quasi--homogeneous, but little more could be said about them. Hence we looked for a more ellegant and concise way of characterizing torsion structures in these cases.

On this line, we tried to gather together this kind of result with the so--called Tate normal form, which was already used by us and M.A. Olalla \cite{GOT} to give a highly efficient algorithm for computing $T(E(\QQ)$ and also by Bennet and Ingram \cite{BI,I} for pointing out rather surprising results concerning $T(E(\QQ)$ related to the pair $(A,B)$.

Let us recall some basic facts about Tate normal form which will be helpful in the sequel. The main result concerning Tate normal form can be stated as follows \cite{Husemoller,Kubert}:

\vs

\noindent {\bf Theorem.--} Every elliptic curve over the rationals with a point $P$ of order $n=4,...,10,12$ can be written in Tate normal form
$$
Y^2 + (1-c) XY  - b Y = X^3 - bX^2,
$$
with the following relations:
\begin{itemize}
\item[(1)] If $n=4$, $b = \alpha, \; c =0$.
\item[(2)] If $n=5$, $b = \alpha, \; c = \alpha$.
\item[(3)] If $n=6$, $b = \alpha + \alpha^2, \; c = \alpha$.
\item[(4)] If $n=7$, $b = \alpha^3 - \alpha^2, \; c = \alpha^2
- \alpha$.
\item[(5)] If $n=8$, $b = (2\alpha - 1)(\alpha - 1), \;
c= b/\alpha$.
\item[(6)] If $n=9$, $c=\alpha^2 (\alpha-1), \; b = c(\alpha
(\alpha-1)+1)$.
\item[(7)] If $n=10$, $c= (2\alpha^3-3\alpha^2+\alpha)/\left[
\alpha - (\alpha -1)^2 \right], \; b=c\alpha^2/\left[ \alpha-
(\alpha-1)^2 \right]$.
\item[(8)] If $n=12$, $c= (3\alpha^2-3\alpha+1)(\alpha-
2\alpha^2)/ (\alpha-1)^3, \; b= c(2\alpha-2\alpha^2-1)/(\alpha-1)$.
\end{itemize}

\vs

\noindent {\bf Remark.--} Basically, our algorithm \cite{GOT} goes as follows: if one wants to know whether or not there is a point of order $n$ on $E(\QQ)$ one considers the corresponding Tate normal form
$$
Y^2 + \left( 1-c_n(\alpha) \right) - b_n(\alpha)Y = X^3 - b_n(\alpha) X^2,
$$
and takes it to a Weierstrass short form
$$
Y^2 = X^3 + A_n (\alpha) X + B_n (\alpha), \mbox{ with } A_n,B_n \in Q(\alpha).
$$

Then our curve has a point of order $n$ if and only if there exists a pair $(u,\alpha) \in \QQ^2$ verifying
$$
u^4 A = A_n (\alpha), \quad u^6 B = B_n (\alpha).
$$

These $A_n,B_n \in \QQ(\alpha)$ become more and more complicated, as $n$ gets bigger; but can be computed once and for all. Hence, the existence of a rational point of order $n$ can be also read in terms of the existence of a solution of the system
$$
\left\{ \begin{array}{rcl}
u^4 A &=& A_n (\alpha) \\
u^6 B &=& B_n (\alpha)
\end{array} \right.
$$

Our concern with this kind of systems was that the solutions had to be searched for in $\QQ$, instead of $\ZZ$. Hence we tried to work with these systems in order to get new ones, where only integral parameters had to be considered. 

\vs

Following this target, our main result in this paper is the following:

\vs

\noindent {\bf Theorem.-- } Given an elliptic curve $E:Y^2=X^3+AX+B$ with $A,B \in \ZZ$, for $n \in \{5,7,8,9\}$ there are homogeneous binary polynomials $F_n,G_n \in \ZZ[p,q]$, at least one of them irreducible, such that $E$ has a rational point of order $n$ if and only if there is a solution to the system
$$
6^4A=F_n, \quad\quad 6^6B=G_n.
$$ 

\vs

\noindent {\bf Remark.--} Rational (actually integral) points of order $n$ have coordinates
$(x_n,y_n)$ which can also be written as homogeneous binary polynomials in the
same variables $\{p,q\}$. 

The following table shows the degrees of $F_n$, $G_n$, $x_n$ and $y_n$ with respect to $\{p,q\}$ for $n \in \{5,7,8,9\}$:

\VS

\begin{center}
\begin{tabular}{ccccc}
$n$ & $\deg(F_n)$ & $\deg(G_n)$ & $\deg(x_n)$ & $\deg(y_n)$ \\
\hline \\
$5$ & $4$ & $6$ & $2$ & $3$ \\
$7$ & $8$ & $12$ & $4$ & $6$ \\
$8$ & $8$ & $12$ & $4$ & $6$ \\
$9$ & $12$ & $18$ & $6$ & $9$ \\
\end{tabular}
\end{center}

\Vs

\noindent {\bf Remark.--} The polynomials $F_n$ and $G_n$ are obtained quickly from $A_n(\alpha)$ and $B_n (\alpha)$, respectively. The important point in the proof is showing that the denominator of $u$ is, at most six (in fact, we show that is either two or three).

\Vs

\noindent {\bf Remark} A useful tool in our proof will be elimination theory for obtaining the explicit expressions of the torsion points. This could as well be obtained from Tate normal form, but dealing a greater amount of dull computation. It has now become customary to regard ellimination theory as an application of Gr\"obner bases, and so we will deal with it. Lots of books cover this issue by now; \cite{Cox} will be a perfect reference for what is needed.

\Vs

\noindent {\bf Remark.--} This result bonds elliptic curves with rational
torsion of order bigger than four to diohantine systems of Thue equations (binary irreducible equations of degree greater than $2$).

These sort of equations have been thoroughly studied along the last
century. To give a brief outline, in 1909, Thue \cite{Thue} proved that 
they have only finitely many solutions. In the sixties, Baker \cite{Baker,DEQ}
gave a theoretic algorithm to find those solutions. Later, in
1989, Tzanakis and Weger \cite{TdW, Smart} combined diophantine
aproximation computational techniques with Baker's theory, to give
a general practical algorithm for solving Thue equations.

The connection with Thue equations will be useful for proving the following corollary:

\vs

\noindent {\bf Corollary.--} Let $n>1$ be an integer. The number of non--isomorphic curves with a point of order $n$ and discriminant $\Delta$ is bounded by an integer depending only on the number of prime factors of $\Delta$.

\Vs

\noindent {\bf Remark.--} Finally, we will mention two problems we are currently working on, for which we hope these results will be helpful:

\begin{enumerate}
\item[(a)] How the torsion subgroup is affected when we consider an algebraic field extension $\QQ \subset \QQ[\alpha]$.
Some results on this line are now known \cite{Kwon,Fujita,Fujita2}, specially for extensions of degree $2^n$ and non--cyclic torsion subgroups, but much is yet to be done.
\item[(b)] The relationship (if any) between the Galois group of $X^3+AX+B$ and the torsion subgroup of $Y^2=X^3+AX+B$. While this relationship is trivial for even order groups, it is not clear at all for the odd case. The given explicit description of $A$ and $B$ may shed some light.
\end{enumerate}

\section{Proof of the main theorem}

We will prove that the system having a solution is a necessary condition, since it will easily checked to be sufficient in each case, by displaying the set of torsion points.

\VS

\noindent \textit{Case} $n=5$. First we will prove that
$$ 
\begin{array}{rcl}
A=    F_5 & = & -27(q^4-12q^3p+14q^2p^2+12p^3q+p^4),\\
B=    G_5 & = & 54(p^2+q^2)(q^4-18q^3p+74q^2p^2+18p^3q+p^4),\\
\end{array}
$$
with $p,q \in \mathbb{Z}\setminus \{0\}$.

Moreover, we will see that the rational points of order five are
$$ 
\left( 3(p^2-6pq+q^2),\pm 108p^2q), \  (3(p^2+6pq+q^2),\pm 108pq^2 \right).
$$

From the proof of the main theorem on \cite{GT}, we
know that if there is a rational point of order five in $E$, then
$$ 
A = -x_1^2-x_1x_2-x_2^2 + (x_1-x_2)x_3,
$$
$$ 
B = \frac{-1}{4}(x_1+x_2)(-3x_1^2+2x_1x_2-3x_2^2+
2(x_1-x_2)x_3).
$$
with 
$$
x_3^2=(2x_1+x_2)(x_1+2x_2)
$$

Besides, the full list of torsion points of order five is
$$
\left\{ \left(x_1, \, \pm \frac{(x_1-x_2)x_4}{2} \right), \; \left( x_2, \, \pm \frac{(x_1-x_2)t}{2} \right) \right\},
$$
where
$$
t^2=3x_1-2x_3+3x_2, \quad\quad x_4^2=3x_1+2x_3+3x_2.
$$

On the other hand, from \cite{GOT} we know that $E$ must be equivalent to
the short Weierstrass form of 
$$
Y^2-\alpha XY-\alpha Y=X^3-\alpha X^2
$$ 
with $\alpha \in \QQ$, that is, 
$$
Y^2=X^3+A_5(\alpha)X+B_5(\alpha),
$$ 
where
$$
A_5(\alpha)=-27-324\alpha -378\alpha ^2+324\alpha ^3-27\alpha ^4,
$$
$$
B_5(\alpha)=54+972\alpha +4050\alpha ^2+4050\alpha ^4-972\alpha ^5+54\alpha ^6.
$$

Hence, there is $u \in \QQ$ such that 
$$ 
A(x_1,x_2,x_3)=u^4 A(\alpha), \quad\quad B(x_1,x_2,x_3)=u^6 B(\alpha),
$$ 
and using Gr\"obner basis we get
$$
x_1=3u^2(\alpha ^2-6\alpha+1),\quad x_2=3u^2(\alpha ^2+6\alpha+1),\quad
x_3=-9u^2(\alpha ^2-1);
$$ 
where $u, \alpha \in \QQ$. 

Now we set
$$
u=\frac{u_1}{u_2},\quad \alpha= \frac{p}{q},
$$ 
with $u_1, u_2, p, q \in \ZZ$, such that $\gcd(p,q)=1=\gcd(u_1,u_2)$, and we will study all
possibilities for $u_1, u_2, p \mbox{ and } q$.

From the last expression of $x_1$, since $\gcd(p^2-6pq+q^2,q^2)=1,$
we have that $q^2$ divides $3u_1^2,$ that is, $q$ divides $u_1$.

If we now substitute those expressions for $x_1,x_2$ and $x_3$ in
$A$ and $B$, and denote $k=u_1/q$, we obtain
$$ 
\begin{array}{rcl}
    A & = & -27(k/u_2)^4(q^4-12q^3p+14q^2p^2+12p^3q+p^4),\\
    B & = & 54(k/u_2)^6(p^2+q^2)(q^4-18q^3p+74q^2p^2+18p^3q+p^4).\\
\end{array}
$$

To finish this proof we just have to show that $u_2=1$ and change $(p,1)$ by $(kp,kq)$. Note that 
$x_1,x_2 \in \ZZ$ becuase of the Nagell--Lutz Theorem \cite{Lutz, Nagell}; and hence
$x_3 \in \ZZ$ as well. From
$$
x_1=\frac{3k^2(p^2-6pq+q^2)}{u_2^2}, \quad x_2=\frac{3k^2(p^2+6pq+q^2)}{u_2^2}, \quad x_3=\frac{9k^2(p^2-q^2)}{u_2^2},
$$
since $\gcd(k,u_2)=1$ and $x_1,x_2,x_3 \in \ZZ$, we have 
$$
u_2^2 \ | \ 3(p^2-6pq+q^2), \quad u_2^2 \ | \ 3(p^2+6pq+q^2),
\quad u_2^2 \ | \ 9(p^2-q^2).
$$ 

Thus $u_2^2$ divides $6(p^2+q^2)$ and $9(p^2-q^2)$. We will consider several
cases:

\begin{enumerate}
\item[(a)] When neither $2$ nor $3$ divides $u_2$, since $u_2^2$ divides $p^2+q^2$
and $p^2-q^2$, we have that $u_2^2$ divides $2p^2$ and $2q^2,$
that is, $u_2=1$.

\item[(b)] If $u_2 \neq 2$ and $2$ divides $u_2$ but $3$ does not, let $u'_2$
be such that $u_2=2u'_2.$ Arguing as above, now with $u'_2$, we
get $u'_2=1$ which is impossible.

\item[(c)] We have an analogous case when $u_2 \neq 3$ and $3$ divides $u_2$
but $2$ does not, setting now $u_2=3u'_2$. And the same when $2$
and $3$ divide $u_2$ and $u_2 \neq 6$, taking $u_2=6u'_2$.

\item[(d)] If $u_2 = 2$, since $B \in \ZZ$ and $\gcd(k,u_2)=1$, it must be
$$
(p^2+q^2)(q^4-18q^3p+74q^2p^2+18p^3q+p^4) \equiv 0 \; \mod \ 32,
$$ 
but this is easily checked to be impossible when $\gcd(p,q)=1$.

\item[(e)] For $u_2 = 3$, since $x_1 \in \ZZ$, it should be 
$$
(p^2-6pq+q^2) \equiv 0 \; \mod\ 3,
$$ 
which is also impossible. And the same occurs if $u_2 = 6$ for 
$$
(p^2-6pq+q^2) \equiv 0 \; \mod\ 12.
$$
\end{enumerate}

Therefore, $u_2=1$ and we get the stated expressions for $A$ and $B$
Once the corresponding substitutions are performed, the rational points of order five are
easily checked to be the ones we stated above.

\VS

\noindent \textit{Case} $n=7$. Now we will show that
$$
\begin{array}{rcl}
 A & = & -27k^4(p^2-pq+q^2)(q^6+5q^5p-10q^4p^2-15q^3p^3+30q^2p^4\\
     &   & \quad -11qp^5+p^6),\\
 B & = & 54k^6(p^{12}-18p^{11}q+117p^{10}q^2-354p^9q^3+570p^8q^4-486p^7q^5\\
     &   & \quad +273p^6q^6-222p^5q^7+174p^4q^8-46p^3q^9-15p^2q^{10}+6pq^{11}+q^{12})\\
 \end{array} 
$$
with $p,q \in \mathbb{Z}\setminus \{0\}$ verifying $p\neq q$, and $k\in \{1,1/3\}$.

Furthermore, we will see that the rational points of order seven are
$$
\begin{array}{l}
    (3k^2(p^4-6p^3q+15p^2q^2-10pq^3+q^4),\pm 108k^3(p-q)^3pq^2),\\
    (3k^2(p^4-6p^3q+3p^2q^2+2pq^3+q^4),\pm 108k^3(p-q)p^2q^3),\\
    (3k^2(p^4+6p^3q-9p^2q^2+2pq^3+q^4),\pm 108k^3(p-q)^2p^3q).
    \end{array}
$$

As we proved in \cite{GT}, if there is a rational point of order seven in
$E$,then 
$$
A=-x_1^2-x_2^2-x_1x_2+ x_4(x_1-x_2), 
$$
$$
4B=(3x_1^3+x_3x_1^2+3x_2^2x_1+x_2^2x_3-2x_1x_2x_3+2x_2^3+ 2(x_2^2-x_1^2)x_4), 
$$
with
$$
x_4^2=(x_2+2x_1)(x_1+x_3+x_2).
$$ 

Again from \cite{GOT}, we know that $E$ must be equivalent to the short Weierstrass form of
the Tate normal form for $n=7$, whose coefficients are
$$
\begin{array}{rcl}
A_7(\alpha) & = & -27(\alpha^2-\alpha+1)(\alpha^6-11\alpha^5+30\alpha^4-15\alpha^3-10\alpha^2+5\alpha+1),\\
B_7(\alpha) & = &
54+324\alpha-810\alpha^2-2484\alpha^3+9396\alpha^4-11988\alpha^5+14742\alpha^6\\
    &  & -26244\alpha^7+30780\alpha^8-19116\alpha^9+6318\alpha^{10}-972\alpha^{11}+54\alpha^{12}.\\
\end{array}
$$

Therefore, there exists $u \in \QQ$ verifying 
$$ 
A(x_1,x_2,x_4)=u^4 A(\alpha),\quad\quad B(x_1,x_2,x_3,x_4)=u^6 B(\alpha).
$$ 

Now, elliminating as above, we get
$$
\begin{array}{rcl}
 x_1 &=& 3u^2(\alpha^4-6\alpha^3+15\alpha^2-10\alpha+1),\\
 x_2 &=& 3u^2(\alpha^4-6\alpha^3+3\alpha^2+2\alpha+1),\\
 x_3 &=& 3u^2(\alpha^4+6\alpha^3-9\alpha^2+2\alpha+1),\\
 x_4 &=& 9u^2(\alpha^2-\alpha+1)(\alpha^2-3\alpha+1).\\
\end{array}
$$

We set again 
$$
u= \frac{u_1}{u_2},\quad \alpha= \frac{p}{q},
$$
with $\gcd(p,q)=1=\gcd(u_1,u_2)$, and we will study all possibilities
for these $u_1, u_2, p, q \in \ZZ.$

First, we must recall from \cite{GT} that $x_1,x_2,x_3 \in \ZZ$, since they are the first coordinates of the six points of order seven on $E(\QQ)$ (Nagell--Lutz again). Then, as
$$
\gcd(p^4-6p^3q+15p^2q^2-10pq^3+q^4,q^4)=1,
$$ 
we have that $q^4\ | \ 3u_1^2$, thus $q^2 \ | \ u_1$. Therefore
$$ 
\begin{array}{rcl}
A & = & -27(u_3/u_2)^4(p^2-pq+q^2)(q^6+5q^5p-10q^4p^2\\
  &   & \quad -15q^3p^3+30q^2p^4-11qp^5+p^6),\\
B & = & 54(u_3/u_2)^6(p^{12}-18p^{11}q+117p^{10}q^2-354p^9q^3\\
  &   & \quad +570p^8q^4-486p^7q^5+273p^6q^6-222p^5q^7\\
  &   & \quad +174p^4q^8-46p^3q^9-15p^2q^{10}+6pq^{11}+q^{12}).\\
\end{array}
$$
where $u_3=u_1/q^2 \in \ZZ$.

Since $\gcd(u_2,u_3)=1$, $u_2^2$ must divide $36(p^2q^2-pq^3)$
because $x_1,x_2 \in \ZZ$, and $u_2^2$ must divide $36(p^2q^2-p^3q)$
because $x_2,x_3 \in \ZZ$. Thus we have that 
$$
u_2^2 \ | \ 36(p^2q-pq^2).
$$ 
But it is also true that
$$
u_2^2 \ | \ 9(p^4-4p^3q+5p^2q^2-4pq^3+q^4),
$$ 
because $x_4 \in \ZZ$. Hence 
$$
u_2^2 \ | \ \gcd(36pq(p-q),9((p-q)(p^3-3p^2q+2pq^2-2q^3)-q^4))=9,
$$ 
that is $u_2= \pm 1$ or $u_2= \pm 3$, so the result is proved. The computation of the points is again direct.

\VS

\noindent \textit{Case} $n=8$. We will show that
$$
\begin{array}{rcl}
A & = & -27k^4(q^8-16pq^7+96p^2q^6-288p^3q^5+480p^4q^4-448p^5q^3+224p^6q^2 \\
    &   & \quad -64p^7q+16p^8) \\
B & = & -54k^6(8p^4-16p^3q+16p^2q^2-8pq^3+q^4)(8p^8-32p^7q-80p^6q^2\\
    &   & \quad +352p^5q^3-456p^4q^4+288p^3q^5-96p^2q^6+16pq^7-q^8)
\end{array}
$$
with $p,q \in \ZZ$ verifying $p \neq q$, $2p \neq q$ and $k \in \{ 1,1/2 \}$.
Given a solution $(p,q)$ of the above system, the points of order eight in the curve are
$$
\left\{ \ \left( 3k^2(-4p^4+20p^3q-20p^2q^2+4pq^3+q^4), \ \pm 108k^3pq(q-p)^3(q-2p) \right),
\right.
$$
$$
\left. \ \left( 3k^2(-4p^4-4p^3q+16p^2q^2-8pq^3+q^4), \ \pm 108k^3p^3q(q-p)(q-2p) \right),
\right\}
$$

As we proved in \cite{GT} the existence of points of order eight impplies the following parametrization:
$$
\begin{array}{rcl}
A(z_1,z_2) &=& \displaystyle \left( -3 z_1^2  + 6 z_1 z_2^2  - 2 z_2^4 \right) \\
B(z_1,z_2) &=& \displaystyle \left( 2z_1-z_2^2 \right) \left( z_1^2+2z_1z_2^2 - z_2^2 \right) \\
0 &=& z_3^2+z_4^2-3z_1 \\
0 &=& z_4^2-z_2(2z_3+z_2)=0
\end{array}
$$
with $z_1,z_2,z_3,z_4 \in \ZZ$.

On the other hand, $E$ must be equivalent to a short Weierstrass form 
$$
Y^2=X^3+A_8(\alpha)X+B_8(\alpha),
$$ 
where
$$
\begin{array}{rcl}
A_8 (\alpha) &=& (-27/\alpha^4)(16\alpha^8-64\alpha^7+224\alpha^6-448\alpha^5 + 480\alpha^4 \\
&& \quad -288\alpha^3+96\alpha^2-16\alpha+1) \\
B_8 (\alpha) &=& (-54/\alpha^6)(64\alpha^{12}-384\alpha^{11}+3520\alpha^9-10296\alpha^8+15840\alpha^7 \\
&& \quad -15568\alpha^6+10272\alpha^5-4560\alpha^4+1328\alpha^3-240\alpha^2+24\alpha-1)
\end{array}
$$

As above, 
$$
A(z_1,z_2) = u^4 A_8 (\alpha), \quad B(z_1,z_2) = u^6 B_8 (\alpha),
$$
yields, after ellimination,
$$
\begin{array}{rcl}
z_1 &=& 3u^2(20\alpha^4-40\alpha^3+28\alpha^2-8\alpha+1)/\alpha^2 \\
z_2 &=& 3u(2\alpha-1)^2/\alpha \\
z_3 &=& 6u(1-\alpha) \\
z_4 &=& 3u(1-2\alpha)/\alpha
\end{array}
$$

We set, as usual $u = u_1/u_2$, $\alpha=p/q$, and perform the corresponding substitution in $z_1$ and $z_2$ to obtain:
$$
\begin{array}{rcl}
A & = & -27(u_1/pqu_2)^4(q^8-16pq^7+96p^2q^6-288p^3q^5+480p^4q^4 \\
  &   & \quad -446p^5q^3+224p^6q^2-64p^7q+16p^8) \\
B & = & -54(u_1/pqu_2)^6(8p^4-16p^3q+16p^2q^2-8pq^3+q^4)(8p^8-32p^7q\\
  &   & \quad -80p^2q^2+352p^5q^3-456p^4q^4+288p^3q^5-96p^2q^6+16pq^7-q^8)
\end{array}
$$

Now, since $z_3 \in \ZZ$ and $\gcd(u_1,u_2)=1$, we have $u_2 \ | \ 6(q-p)$. Also $z_4 \in \ZZ$, therefore $u_2 \ | \ 3(q-2p)$. This proves $u_2$ divides both $6p$ and $3q$, hence $u_2=1$ or $u_2=3$. 

On the other hand, $z_3 \in \ZZ$ and $\gcd(p,q-p)=1$ imply $q\ | \ 6u_1$; and $A \in \ZZ$ imply $p^4\ | \ 27u_1^4$, hence $p \ | \ u_1$.

The equations involving $A$ and $B$ are now proved, but it remains showing $k \in \{1,1/2\}$. From above, we actually need to show that $3$ does not divide 
$$
q^8-16pq^7+96p^2q^6-288p^3q^5+480p^4q^4-446p^5q^3+224p^6q^2-64p^7q+16p^8.
$$

We will consider two cases:

\begin{enumerate}
\item[(a)] If either $3 \ | \ q$ or $3\ | \ p$, then as $\gcd(p,q)=1$ the result is straightforward.
\item[(b)] Otherwise we have 
$$
p \equiv \pm 1 \; \mod \ 3,
$$
and so does $q$. Checking all possibilities it is easy to see that $3$ does not divide the above expression.
\end{enumerate}

\VS

\noindent \textit{Case} $n=9$. Finally we will prove that
$$
\begin{array}{rcl}
A & = &-27k^4(q^3-3p^2q+p^3)(q^9-9q^7p^2+27q^6p^3-45q^5p^4+54q^4p^5\\
    &   & \quad -48q^3p^6+27p^7q^2-9p^8q+p^9), \\
B & = &54k^6(p^{18}-18p^{17}q+135p^{16}q^2-570p^{15}q^3+1557p^{14}q^4-2970p^{13}q^5\\
    &   & \quad +4128p^{12}q^6-4230p^{11}q^7+3240p^{10}q^8-2032p^9q^9+1359p^8q^{10} \\
    &   & \quad -1080p^7q^{11}+735p^6q^{12}-306p^5q^{13}+27p^4q^{14}+42p^3q^{15}\\
    &   & \quad -18p^2q^{16}+q^{18}),\\
\end{array}
$$
with $p,q \in \mathbb{Z}\setminus \{0\}$ verifying $p\neq q$, and $k=1$ or $k=1/3$.

Moreover, we will show that the rational points of order 9 are
$$
\begin{array}{l}
     (3k^2(p^6+6p^5q-15p^4q^2+14p^3q^3-6p^2q^4+q^6),\\
       \ \ \ \ \pm 108k^3p^4q(p^4-3p^3q+4p^2q^2-3pq^3+q^4)),\\
     (3k^2(p^6-6p^5q+21p^4q^2-34p^3q^3+30p^2q^4-12pq^5+q^6),\\
       \ \ \ \ \pm 108k^3pq^2(p^6-5p^5q+11p^4q^2-14p^3q^3+11p^2q^4-5pq^5+q^6)),\\
     (3k^2(p^6-6p^5q+9p^4q^2-10p^3q^3+6p^2q^4+q^6),\\
       \ \ \ \ \pm 108k^3p^2q^4(p^3-2p^2q+2pq^2p-q^3)).\\
\end{array}
$$

In this case, we have 
$$
A=27z_1^4+6z_1z_2,\qquad B=z_2^2-27z_1^6,
$$ 
with $z_1,z_2 \in \ZZ$, and $E$ must be equivalent to the
short Weierstrass form of
$$
Y^2+(1-{\alpha}^2(\alpha-1))XY-{\alpha}^2(\alpha-1)({\alpha}^2-\alpha+1)Y=X^3-{\alpha}^2(\alpha-1)({\alpha}^2-\alpha+1)X^2
$$
that is, to $Y^2=X^3+A(\alpha)X+B(\alpha),$ where
$$
\begin{array}{rcl}
A(\alpha) & = & -27\alpha^{12}+324\alpha^{11}-1458\alpha^{10}+3456\alpha^9-5103\alpha^8\\
            &   & \quad +4860\alpha^7-3078\alpha^6+972\alpha^5+486\alpha^4-756\alpha^3\\
            &   & \quad +324\alpha^2-27,\\
B(\alpha) & = & 54\alpha^{18}-972\alpha^{17}+7290\alpha^{16}-30780\alpha^{15}+84078\alpha^{14}\\
            &   & \quad -160380\alpha^{13}+222912\alpha^{12}-228420\alpha^{11}+174960\alpha^{10}\\
            &   & \quad -109728\alpha^9+73386\alpha^8-58320\alpha^7+39690\alpha^6\\
            &   & \quad -16524\alpha^5+1458\alpha^4+2268\alpha^3-972\alpha^2+54,\\
\end{array}$$
with $\alpha \in \QQ$. Hence, if the existence of a rational point
$P$ of order 9 implies to the existence of $u, \alpha \in \QQ$ such that
$$ 
\frac{A(z_1,z_2)}{A(\alpha)} = u^4,\ \ \frac{B(z_1,z_2)}{B(\alpha)} =u^6,
$$
and this occurs if and only if
$$ 
\begin{array}{rcl}
     z_1 & = & u(1-3{\alpha}^2+{\alpha}^3), \\
     z_2 & = & -9z_1^3+108u^3{\alpha}^3(\alpha-1)^3,\\
 \end{array}
 $$
with $u,\ \alpha \in \mathbb{Q}.$

Recall from \cite{GT} there exists a polynomial $R_9 \in \ZZ[z_1,z_2,z_3]$ of degree nine in $z_3$ with 
three integer roots (w.r.t. $z_3$, for a given curve with order nine torsion), which were the first 
coordinates of the rational points or order 9 in the curve, that is, $x(P),\ x(2P)$ and $x(3P)$. If we
substitute this expressions for $z_1$ and $z_2$ in $R_9$ we
get only the following rational (hence integral) roots:
$$
\begin{array}{rcl}
    x(P) & = & 3z_1^2-4(3u)^2{\alpha}^2(\alpha-1), \\
    x(2P) & = & 3z_1^2+4(3u)^2{\alpha}^3(\alpha-1)^2,\\
    x(4P) & = & 3z_1^2+4(3u)^2\alpha(\alpha-1)^3.\\
\end{array} 
$$

Again we set $u=u_1/u_2$, $\alpha=p/q,$ with $\gcd(p,q)=1=\gcd(u_1,u_2),$ and will study all
possibilities for these integers.

From 
$$
x(P)=3\frac{u_1^2}{u_2^2} \left( \frac{p^6-6p^5q+9p^4q^2-10p^3q^3+6p^2q^4+q^6}{q^6} \right) \in
\ZZ
$$ 
we get $u_3=u_1/q^3 \in \ZZ$, because
$$
\gcd(p^6-6p^5q+9p^4q^2-10p^3q^3+6p^2q^4+q^6,q^6)=1.
$$ 
Now, since $x(P), z_1 \in \ZZ$,we have
$$
x(P)-3z_1^2=-2^2{\cdot}3^2u_3^2q^3\frac{p^2(p-q)}{u_2^2}\in \ZZ.
$$

If we set $u_2=2^a3^bu_4$, with $a=1$ if $2$ divides $u_2$ and
$a=0$ if not, and the same for $b$ and $3$, we get
$$ 
\frac{-2^2{\cdot}3^2p^2(p-q)}{2^{2a}3^{2b}u_4^2} \in \ZZ.
$$

There are several cases:
\begin{itemize}
\item[(a)] If $a=b=0$ then $u_2^2$ divides $p^2(p-q)$. Since
$$
p^6-6p^5q+9p^4q^2-10p^3q^3+6p^2q^4+q^6 \equiv 0 (\mod\ u_2^2),
$$
we have 
$$
p^2(p-q)(p^3-5p^2q+4pq^2-6q^3)\equiv -q^6 (\mod\ u_2^2),
$$
therefore $q^6 \equiv 0 (\mod\ u_2^2)$. But then, since $q^3$
divides $u_1$, we get $u_2=1$.

\item[(b)] If $a=1, b=0$ and $u_2 \neq 2$ then $u_2=2u_4$ with
$u_4\neq 1$, and now $u_4^2$ divides $p^2(p-q)$. As
above we get $u_4=1$ which is impossible.

\item[(c)] If $a=0, b=1$ and $u_2 \neq 3$, then $u_2=3u_4$ with
$u_4\neq 1$ and $u_4^2$ divides $p^2(p-q)$. Again we get a
contradiction because $u_4$ divides $q^3$.

\item[(d)] If $a=1, b=1$ and $u_2 \neq 6$, then $u_2=6u_4$ with
$u_4\neq 1$ and it follows the same as in former cases.

\item[(e)] If $u_2=2$ we get 
$$
p^6-6p^5q+9p^4q^2-10p^3q^3+6p^2q^4+q^6\equiv 0 (\mod\ 2^2),
$$ 
but there exist no $p,q \in \ZZ$ relatively primes satisfying that congruence.

\item[(f)] If $u_2=6$ we get 
$$
p^6-6p^5q+9p^4q^2-10p^3q^3+6p^2q^4+q^6 \equiv 0 (\mod\ 12),
$$ 
and again there is no solution.
\end{itemize}

Hence we have proved that $u_2=1$ or $u_2=3$ as wanted. Moreover,
if we substitute $A$ and $B$ by its expressions with $k$, $p$ and
$q$ in our polynomial of degree 9, we get three integer roots,
namely:
$$\begin{array}{rcl}
x(P) &=& 3k^2(p^6+6p^5q-15p^4q^2+14p^3q^3-6p^2q^4+q^6),\\
x(2P)&=&3k^2(p^6-6p^5q+21p^4q^2-34p^3q^3+30p^2q^4-12pq^5+q^6),\\
x(4P)&=&3k^2(p^6-6p^5q+9p^4q^2-10p^3q^3+6p^2q^4+q^6).\\
\end{array}$$

\VS

\noindent {\bf Remark.--} One may ask whether $k$ is actually necessary for $n=7,8,9$. The following curves show it is. For all of these curves, the system $\{ A=F_n, \; B= G_n\}$ does not have an integral solution, while $\{ 6^4A=F_n, \; 6^6 B= G_n\}$ has.

\Vs

\begin{tabular}{cl}
{\bf Torsion order} & $\qquad$ {\bf Curve} \\
\hline \\
$7$ & $Y^2 = X^3-43X+166$ \\
$8$ & $Y^2 = X^3-22187952X+23592760704$ \\
$9$ & $Y^2 = X^3-219X+1654$
\end{tabular}

\section{Proof of the corollary}

We will prove now the corollary stated in the introduction: For a given $n$, a prime power, the number of non--isomorphic curves with a point of order $n$ and discriminant $\Delta$ is bounded by an integer depending only on the number of prime factors of $\Delta$. 

For every such $n$, we have an explicit expression of $A$ and $B$ in terms of two variables, either from \cite{GT} or from the above theorem, hence we can compute $\Delta$. Here we list the results (variables are now $x$ and $y$):

\Vs

\begin{tabular}{clc}
$n$ & $\Delta$ & {\bf Degree} \\
\hline \\
$2$ & $2^4(4x-y^2)(x+2y^2)^2$ & $3$ \\
$3$ & $2^43^3(5x^3+y)(9x^3+y)^3$ & $4$ \\
$4$ & $2^4y^2(12x-5y^2)(3x-y^2)^4$ & $6$ \\
$5$ & $2^{12}3^{12}x^5y^5(x^2+11xy-y^2)$ & $12$ \\
$7$ & $2^{12}x^7y^7(x^3-8x^2y+5xy^2+y^3)(y-x)^7$ & $24$ \\
$8$ & $3^{12}x^8y^2(8x^2-8xy+y^2)(2x-y)^4(x-y)^8$ & $24$ \\
$9$ & $-2^8x^9(x^3-6x^2y+3xy^2+y^3)(x^2-xy+y^2)^3(x-y)^9$ & $27$ 
\end{tabular}

\Vs

\noindent where Degree stands for the homogeneous degree obtained by doing
\begin{eqnarray*}
n=2 && y^2 \longmapsto y \\
n=3 && x^3 \longmapsto x \\
n=4 && y^2 \longmapsto y
\end{eqnarray*}

After these substitutions we get, for all $n$, Thue equations 
$$
F(x,y) = \Delta
$$
whose set of solutions contains the solutions of the original equations shown above. Now we can apply the Evertse bound \cite{Evertse}: the number of primitive solutions of a Thue equation of degree $r$ given by
$$
F(x,y) = h
$$
is bounded by
$$
7^{15 {\tiny \left[ \left( \begin{array}{c} r \\ 3 \end{array} \right) +1 \right]^2}} + 6 \cdot 
7^{2 {\tiny \left( \begin{array}{c} r \\ 3 \end{array} \right)} (t+1)},
$$
where $t$ is the number of prime factors of $h$. Note that isomorphic curves give rise to scaled solutions, hence by considering primitive solutions we are sure to include all non--isomorphic cases. 

When we apply this bound to our equations, we note that the bound obtained for $n$ prime is smaller than the attached by any multiple of $n$. Hence we finally get the following result:

\Vs

\noindent {\bf Corollary.--} Let $n>1$ be an integer. The number of non--isomorphic curves with a point of order $n$ and discriminant $\Delta$ is bounded by an integer $M_n (t)$ depending only on $t$, the number of prime factors of $\Delta$. More precisely:

\begin{eqnarray*}
n=2,4,6,8,10,12 &&  M_n(t) = M_2(t) = 7^{60} + 6\cdot 7^{2(t+1)} \\
n=3,9 && M_3(t) = M_9(t) = 7^{375} + 6\cdot 7^{8(t+1)} \\
n=5 && M_5(t) = 7^{1815} + 6 \cdot 7^{20(t+1)} \\
n=7 && M_7(t) = 7^{19440} + 6 \cdot 7^{70(t+1)}
\end{eqnarray*}

\end{document}